\input amstex
\documentstyle{amsppt}
\magnification=\magstep1
\vsize 21 true cm
\hsize 15.7 true cm
\abovedisplayskip=1in
\belowdisplayskip=1in
\hoffset=0.15in
\voffset=0.05in

\catcode `\@=11
\let \logo@=\relax
\catcode  `\@=\active
\def\deg{\operatorname{deg}}%
\def\dim{\operatorname{dim}}%
\def\reg{\operatorname{reg}}%
\def\codim{\operatorname{codim}}%
\def\Ker{\operatorname{Ker}}%
\def\rank{\operatorname{rank}}%
\def\length{\operatorname{length}}%
\def\Supp{\operatorname{Supp}}%

\parindent=5mm

\topmatter
\title
{ Generic Projection Methods and Castelnuovo Regularity
of Projective Varieties}
\endtitle
\leftheadtext\nofrills{\hskip10.5cm SIJONG  KWAK}
\rightheadtext\nofrills{\hskip-2.5cm GENERIC PROJECTION METHODS AND CASTELNUOVO REGULARITY}
\author
Sijong Kwak
\endauthor

\address
 School of Mathematics, Korea Institute for Advanced Study,
207-43 Chungryangri-dong, Dongdaemoon-gu, Seoul 130-010, Korea
\endaddress

\email
sjkwak\@kias.kaist.ac.kr
\endemail


\abstract
For a reduced, irreducible projective variety $X$ of degree $d$ and
codimension $e$ in $\Bbb P^N$ the Castelnuovo-Mumford regularity
$\reg{X}$ is defined as the least $k$ such that $X$ is $k$-regular,
i.e\. $H^i(\Bbb P^N,\Cal I_X(k-i))= 0$ for $i\ge1$, where $\Cal
I_X\subset\Cal O_{\Bbb P^N}$ is the sheaf of ideals of $X$. There is a
long standing conjecture about $k$-regularity (see \cite{EG}):
$\reg{X}\le d-e+1$. Generic projection methods proved to be
effective for the study of regularity of smooth projective varieties
of dimension at most four (cf\. \cite{BM}, \cite{K2}, \cite{L},
\cite{Pi}, and \cite{R1}) because there are nice vanishing theorems for
cohomology of vector bundles (e.g. the Kodaira-Kawamata-Viehweg
vanishing theorem) and detailed information about the fibers of
generic projections from $X$ to a hypersurface of the same dimension.

Here we show by using methods similar to those used in \cite{K2} that
$\reg{X}\le(d-e+1)+10$ for any smooth fivefold and $\reg{X}\le(d-e+1)+20$
for any smooth sixfold.

Furthermore, using similar methods we give a bound for the regularity
of an arbitrary (not necessarily locally Cohen-Macaulay) projective
surface $X$ of degree $d$ and codimension $e=N-2$ in $\Bbb P^N$. To wit, we show that
$\reg{X}\le(d-e+1)d-(2e+1)$. This is the first bound for surfaces
which does not depend on smoothness.
\endabstract

\affil
School  of  Mathematics \\ Korea Institute for Advanced Study
\endaffil

\endtopmatter

\document

\baselineskip=12.7pt

\head
\bf \S0. Introduction
\endhead

For a given nondegenerate projective variety $X$ of dimension
$n$, codimension $e$ and degree $d$ in $\Bbb P^N$, one can
easily show that $X$ is set-theoretically an intersection of
hypersurfaces (actually cones) of degree at most $d$. Furthermore,
If $X$ is smooth, then $X$ is
scheme-theoretically cut out by homogeneous polynomials of degree
$d$, i.e. there is a surjection ${\bigoplus}_{i=1}^m \Cal
O_{\Bbb P^N}(-d)\rightarrow\Cal I_X$ with $m \ge e$, where $\Cal
I_X\subset\Cal O_{\Bbb P^N}$ is the sheaf of ideals of $X$ \cite{Mu2}.
Then, it is natural to ask whether the degrees of all minimal
generators of the saturated ideal of $X$ are also bounded by $d$.
More strongly, it has been conjectured that the degrees of all minimal
generators are bounded by $d-e+1$. An important role in the study of
this question is played by the Castelnuovo-Mumford regularity
$\reg{X}$. According to \cite{EG}, \cite{Mu1}, $X$ is $m$-regular
iff one of the following conditions holds:
\roster
\item
$H^i(\Bbb P^N ,\Cal I_{X}(m-i))=0$ for all $i\ge1$;
\item
$H^i ({\Bbb P}^N ,  {\Cal  I}_{X}(j)) = 0$ for $i\ge1$, $i+j \ge m$;
\item
For all $k \ge 0$ the degrees of minimal generators of the $k$-th
syzygy modules of the homogeneous saturated ideal ${I}_X$ of $X$ are bounded
by $k+m$.
\endroster

More generally, a  coherent  sheaf  $\Cal F$
on ${\Bbb  P}^N$  is said to be $m$-regular  if
$ H^i (\Bbb P^N ,  \Cal F(m-i))=0$ for all $i>0$, and the regularity
of $\Cal F$ is defined by the formula
$$
\reg{\Cal F}=\min{\{m\in\Bbb Z\:\Cal F\,\,  \text{is}  \,\,
m\text{-regular}\}}.
$$

In particular, $\reg{X}$ is defined as $\reg{\Cal I_X}$. In
general, $\reg{\Cal F}$ may be negative; however, it is not hard to show that
$\reg{X}\ge2$ and $X$ is $2$-regular if and only if $X$ is of minimal degree.

A well known conjecture due to Eisenbud and Goto (cf\. \cite{EG}),
gives a bound for regularity in terms of the $d$ and $e$:
$$
\reg{X} \le d-e+1
$$
A useful tool for the study of regularity of smooth projective
varieties of small dimension is provided by generic projection methods
(cf\. \cite{BM}, \cite{K2}, \cite{L}, \cite{Pi}, and \cite{R1}).
Application of these methods depends on the existence of nice
vanishing theorems for cohomology of vector bundles (e.g. the
Kodaira-Kawamata-Viehweg vanishing theorem) and detailed information
about the fibers of generic projections from $X$ to a hypersurface of
the same dimension. There are good bounds for regularity of smooth
projective varieties of $\dim{X} \le 4$. More precisely, $\reg{X}
\le d-e+1$ for integral curves and smooth surfaces (see \cite{GLP},
\cite {L}) and the author proved the bound $\reg{X}\le(d-e+1)+1$
for smooth threefolds and $\reg{X} \le (d-e+1)+4$ for smooth
fourfolds \cite{K2}. The best known bound for the regularity of a
{\it smooth} projective variety $X$ of {\it arbitrary} dimension $n$ is much worse
than expected, viz\. $\reg{X} \le\min{\{e,n\}}\cdot d-n+1$ (cf\.
\cite{BEL}).

The goal of the present note is to prove new regularity bounds for
arbitrary projective surfaces (Theorem 3.3), for smooth fivefolds (Proposition 2.4, Theorem 2.10) and for smooth sixfolds (Proposition 2.8, Theorem 2.13).
Main idea is to choose special subspaces of the vector spaces of homogenious
polynomials of degree $n\le 7$ in order to show $n$-normality of finite schemes appearing as fibers
of generic projection from a given variety to a hypersurface.


{\bf  Acknowledgements.}
This paper owes a great deal to Greenberg's unpublished Ph.D. thesis
(\cite{G}) written under the direction of Professor H. Pinkham and of
course to Lazarsfeld's important paper(\cite{L}). It is a pleasure to
thank H.~Pinkham  for information on J.~Mather's theorem and
F.~L.~Zak for valuable discussions and his help in improving the clarity
of exposition.

\head
\bf \S1. Basic background
\endhead

In  this  section  we  recall  the
definitions and basic results  which
will be used in subsequent sections. We work over an algebraically
closed field of characteristic zero.

\proclaim{Lemma 1.1}
Let  $X$  be  a  nondegenerate integral scheme of
dimension $n$  in  $\Bbb P^N$,  and let  $Y=X \cap H$  be
a generic hyperplane section.  Then
\roster
\item"\rom{(a)}"
$\reg{Y}\le\reg{X}$\rom{;}
\item"\rom{(b)}"
If $Y$ is $m$-regular, then $X$ is $m+h^1({\Cal I}_{X}(m-1))$-regular.
\endroster
\endproclaim
\demo{Proof}
$(a)$ can be easily checked and $(b)$ is proved in \cite {Mu 2}, page 102.
\qed
\enddemo
Let  $X$ be a nondegenerate zero-dimensional subscheme of length $d$,
not necessarily reduced, and let $N=\dim{\langle X\rangle}$, where
$\langle X\rangle=\Bbb P^N$ is the span of $X$.
Let $\varphi_X$ be the Hilbert function, and let
${\Cal  P}_X$ be the Hilbert polynomial of $X$. It is easy to verify
that the following conditions are equivalent:
\roster
\item"\rom{(i)}" $X$ is $(m+1)$-regular\rom;
\item"\rom{(ii)}" $X$ is $k$-normal for all $k\geq m$\rom;
\item"\rom{(iii)}" $\varphi_X(t)=\Cal P_X(t)=d$ for all $t\geq m$.
\endroster

Let's put
$$
t=\max{\{k\mid\dim{\langle X'\rangle}=\length{X'}-1\quad\forall X'
\subset X,\,\length{X'}\leq k+1\}}.
$$
It is clear that $1\leq t\leq N$, and that $t=1$ iff $X$ has a
trisecant line.

The following Proposition 1.2. and Corollary 1.3. were communicated to me
by F.~L.~Zak. However, for lack of suitable references we give brief
proofs here.

\proclaim{Proposition 1.2} In the above situation,
$X$ is $k$-normal for all $k\ge\lceil \frac{d-N-1}t \rceil+1$, where
$\lceil a \rceil$
is the smallest integer that is not less than $a$.
\endproclaim
\demo{Proof}
We proceed with induction on $N$. If $N=t$,\,\,i.e. $X$ is a
``general position scheme''
then it is proved in Theorem 28.8, \cite {Pe}.
Let's fix an integer $N_0$ and suppose Proposition 1.2 holds
for $t\le N\le N_0-1$. For $N=N_0$, we may also assume that Proposition
1.2 is true for finite schemes of smaller degree than $d$. Let $A$ be a
graded homogeneous ring of $X$. Equivalently, we show the surjectivity of
the natural morphism
$$
A_r\rightarrow H^0(X, \Cal O_X(r))
$$
for all $r$ such that $d\le tr+(N-t)+1$.
Choose a hyperplain $H$ such that $\deg{(X\cap H)}\ge N$ and
$\langle Y \rangle=H$. Consider the following diagram as shown in Theorem 28.8
\cite {Pe}:
$$
\CD
0  @.  0  @.  @.\\
@VVV  @VVV  @.  @.\\
{\lbrack A/(0:H) \rbrack}_r @> \alpha_r >>
H^0(Z, \Cal O_Z(r))   \simeq \Bbb C^{d_1}\\
\hskip-6mm\times H  @VVV  @VVV  @.  @.\\
A_{r+1}
@> \rho_{r+1} >>
H^0(X, \Cal O_X(r+1))
\simeq \Bbb C^d\\
@VVV  @VVV  @.  @.\\
{\lbrack A/HA \rbrack}_{r+1}
@> \beta_{r+1} >>
H^0(Y, \Cal O_Y(r+1))\simeq \Bbb C^{d_2}\\
@VVV  @VVV  @.  @.\\
0  @.  0  @.  @.\\
\endCD
\tag*
$$
where $Y=X\cap H$, $\deg{Y}=d_2$ and $Z$ is the subscheme of $X$
with degree $d_1\ge 1$ corresponding to the graded ring $A/(0:H)$. Clearly,
any closed subscheme of degree $t+1$ in either $Y$ or $Z$ spans $\Bbb P^t$.
So by induction hypothesis,  $\alpha_r$ is surjective for all $r$ such that
$d_1\le tr+(n-t)+1$, $n=\dim{\langle Z \rangle}$ (Note that if $n < t$, then
$Z$ is a ``general position scheme'' and $\alpha_r$ is surjective for all
$r\ge 1$) and $\beta_{r+1}$ is surjective for all $r$ such that
$d_2\le t(r+1)+(N-1-t)+1$. It is easily checked that $d=d_1+d_2\le t(r+1)+(N-t)+1$
implies $d_2\le t(r+1)+(N-1-t)+1$ and $d_1\le tr+(n-t)+1$,
which means $\alpha_r$ and $\beta_{r+1}$ are surjective. By snake lemma,
$\rho_{r+1}$ is also surjective.
\qed
\enddemo

\proclaim{Corollary 1.3} In the same situation as in Proposition 1.2,
\roster
\item"\rom{(a)}" $X$ is $(d-N)$-normal but fails to be $(d-N-1)$-normal if and only if
$X$ has a $(d-N+1)$-secant line\rom{;}
\item"\rom{(b)}" If $t=N$, i.e\. $X$ is a ``general position'' scheme,
then $X$ is $k$-normal for all $k\ge\lceil\frac{d-1}N\rceil$.
\endroster
\endproclaim
\demo{Proof}
For (a), if $X$ has a $(d-N+1)$-secant line then clearly it fails to be
$(d-N-1)$-normal. Conversely, suppose $X$ is $(d-N)$-normal but fails to be
$(d-N-1)$-normal. Then, we proceed with induction on $N$. it is clear for $N=1$.
Suppose that (a) is true for $\dim{\langle X \rangle}\ < N$ and $X$ has no
$(d-N+1)$-secant line. Choose a hyperplain $H$ such that
$N\le \deg{(Y)}\le d-1, \,\,Y=X\cap H$ and $\langle Y \rangle=H$. Since $Y$ has
also no $(d-N+1)$-secant line, $Y$ is $(d-N-1)$-normal by induction hypothesis.
Similarly,
since $1\le \deg{Z}\le d-N$ and we can choose $H$ such that $Z$ has no
$(d-N)$-secant line, $Z$ is $(d-N-2)$-normal. From the same commutative diagram as $(*)$,
we know that $X$ is $(d-N-1)$-normal which contradicts our assumption.
(b) is clear with $t=N$.
\qed
\enddemo

\definition{Definition 1.4} A scheme $X$ is called {\it punctual\/} if
$\operatorname{Supp}{X}=x$, where $x\in X$ is a point. A punctual
scheme $X$ is called {\it curvilinear\/} if ${\Cal O}_x$ is isomorphic
to ${\Bbb C}[x]/(x^k)$ for some $k \ge 1$.
\enddefinition

It is clear that a punctual scheme is curvilinear if and only if it
admits an embedding into a smooth curve.

\proclaim{Lemma 1.5}  Let $X$ be a $n$-dimensional smooth projective
variety in ${\Bbb P}^N$, and suppose that $n=\dim{X}\le5$.  Let
$\Lambda^{N-n-2}$ be a general linear subspace of dimension $(N-n-2)$,
so that, in particular, $\Lambda$ is disjoint from $X$, and let
$\pi_{\Lambda}$ be the projection with center $\Lambda$, and put
$Y=\pi_{\Lambda}(X)\subset\Bbb P^{n+1}$. Then all fibers of
${\pi}_{\Lambda}\:X\to Y$ are curvilinear.
\endproclaim

\demo{Proof}
Let $W \subset V$ be two linear spaces with
$\dim{W}=n$ and $\dim{V}=N$. Then, by an easy computation, the
Schubert cell $\{L \in G(k, N); \dim{(L \cap W)}\geq t\}$ has
codimension $t(N-k-n+t)$ in $G(k, N)$. Therefore, for a given
nonsingular variety $X$ in $\Bbb P^N$, $X_q=\{x \in X:
\dim{(\overline{T_x(X)}}\cap\Lambda)\geq q-1, \Lambda\cap X=
\varnothing, \dim{\Lambda}=N-n-2 \}$  has
codimension $q(q+1)$. Therefore, if $\dim{X}\le5$, then
$X_q=\varnothing$ for $q \ge 2$. This implies that for a generic
projection $\pi_{\Lambda}: X \rightarrow Y \subset \Bbb P^{n+1}$,
$n \le 5$,  $X_q$ is empty for $q\geq2$; in other words, it has only
curvilinear fibers.
\qed
\enddemo

\proclaim{Theorem 1.6 {\rom(J. Mather)}}
Let  $X\subset\Bbb P^N$ be a smooth nondegenerate $n$-dimensional
variety, let $\Lambda^{N-n-2}\subset\Bbb P^N$ be a generic linear
subspace, and let ${\pi}_{\Lambda}\:\Bbb P^N\dashrightarrow\Bbb P^{n+1}$,
$Y=\pi_{\Lambda}(X)\subset\Bbb P^{n+1}$.  Let $Y_k=\{y\in
Y\mid\length{{\pi}_{\Lambda}^{-1}(y)}\geq k\}$, and put
$X_k=\pi^{-1}_{\Lambda}(Y_k)$, so that $X_1\supset\cdots\supset X_k
\supset X_{k+1}\cdots$ is a decreasing filtration. Assume that $n \le
14$, so that we are in Mather's ``nice'' range. Then
$X_{n+2}=\varnothing$ and $\dim{X_{k}}\leq n+1 - k$.  If
$\dim{X_k}=\dim{Y_k}=n+1 - k$, then there exists a dense open subset
of $Y_k$ over which all the fibers of $\pi_{\Lambda}$ are reduced.
\endproclaim

\demo{Proof}
  This follows from the main theorem of \cite{Ma1} and the discussion in
  \S 5 of \cite{Ma2}.
A key ingredient is the inequality
$$
\sum_{x \in \pi^{-1}(y)} (\delta_x + \gamma_x) \le n+1,\qquad y\in Y,
\tag1.0
$$
where $\delta_x=\length{\Cal O_{\pi^{-1}(y),x}}$ and $\gamma_x$ is
another non-negative invariant introduced by J.~Mather for all stable
germs in the ``nice'' range (cf\. \cite{Ma2}); in particular,
$\gamma_x=k-1$  if  ${\Cal O}_x \simeq {\Bbb C}[x]/(x^k)$ for some $k
\ge 1$), which is always the case for $n\leq5$.
\qed
\enddemo

\remark{\bf Remark 1.7} Let $X\subset\Bbb P^N$ be a smooth
nondegenerate $n$-dimensional subvariety and let $S_k(X)$ be the locus
of $k$-secant lines of $X$ in $\Bbb P^N$. Assume that $n\le14$. Then
by Theorem~1.6 one has $\dim{S_{n+2-k}}\le n+1+k$, which gives us some
information on ``collinear'' fibers of a generic linear projection of
$X$ to a hypersurface.
\endremark

\head
\bf \S2. Castelnuovo regularity for smooth varieties of dimension $5$ and $6$
\endhead

Let $X$ be a $n$-dimensional smooth projective variety of degree $d$
and codimension $e$ in $\Bbb P^N$ defined over the  field $\Bbb C$ of
complex numbers. We will use the
general construction considered in \cite{L}, \cite{G}, and \cite{K2}.
Let $\Lambda=\Bbb P^{N-n-2}\subset\Bbb P^N$, $\Lambda\cap
X=\varnothing$, $\Lambda=\Bbb P(V)$ be a general linear subspace, and
let $\pi_{\Lambda}\:X\to Y$ be the projection with center at
$\Lambda$, so that $Y\subset\Bbb P^{n+1}$ is a hypersurface.
Let $\Cal V$ be a collection of linear subspaces $V_j\subset S^j(V)$
such that $V_1 = V$ and $V_2 = S^2(V)$. Consider the natural
restriction morphism $\tilde{\omega}_{n,k,\Cal V}$. If
$\tilde{\omega}_{n,k,\Cal V}$ is surjective, then we get the following
exact sequence:
$$
0\to E_{n,k,\Cal V}\to V_k\otimes\Cal O_{\Bbb P^{n+1}}(-k)\oplus\cdots
\oplus V_1\otimes\Cal O_{\Bbb P^{n+1}}(-1)\oplus\Cal O_{\Bbb P^{n+1}}
\overset{\tilde\omega_{n,k,\Cal V}}\to\longrightarrow
\pi_{\Lambda\,*}\Cal O_X \to0.
\tag2.0
$$

\proclaim{Lemma 2.1} Suppose that $\tilde\omega_{n,k,\Cal V}$ is
surjective. Then
\roster
\item
$\reg{(E_{n,k,{\Cal V}}^*)}\le-2$\rom;
\item
$E_{n,k,{\Cal V}}^*$ is $(-3)$-regular if and only if $X$ is linearly normal
and $H^0(\Cal I_{X/\Bbb P^N}(2))=H^1(\Cal O_X)=0$.
\endroster
\endproclaim

\demo{Proof} (1) is proved in \cite{L, Lemma~2.1}.

(2)  Similar arguments are given in \cite {Al}. By definition,
$E_{n,k,\Cal V}^*$ is $(-3)$-regular iff $H^i(\Bbb P^{n+1}, E_{n,k,\Cal V}^*(-3-i))=0$
for $i>0$. By Serre's duality, this is equivalent to
$H^j(\Bbb P^{n+1}, E_{n,k,\Cal V}(2-j))=0$ for $j\le n$.
For $j\ge3$ this vanishing follows from the Kodaira Theorem.
From the exact cohomology sequence corresponding to~(2.0) it follows that:
$$
\gathered
H^0(\Bbb P^{n+1}, E_{n,k,\Cal V}(2))=0\text{\ if and only if\ }
H^0(\Cal I_{X/\Bbb P^N}(2))=0,\\
H^1(\Bbb P^{n+1}, E_{n,k,\Cal V}(1))=0\text{\ if and only if\
}X\text{\ is linearly normal},\\
H^2(\Bbb P^{n+1},E_{n,k,{\Cal V}})=0\text{\ if and only if\ }
H^1(\Cal O_X)=0.
\endgathered
$$
This completes the proof of (2).

\qed
\enddemo

\remark{\bf Remark 2.2} For varieties of small codimension,
conditions of Lemma~2.1 can be verified using Zak's Linear
Normality Theorem ($X$ is linearly normal if $N<\frac32n+2$)
and Barth's Lefschetz Theorem ($H^1(\Cal O_X)=0$ if
$N<2n$). To verify that $X$ is not contained in a quadric, it suffices
to show that the trisecant lines of $X$ fill up the ambient linear
space $\Bbb P^N$.
\endremark

\proclaim{Lemma 2.3} Suppose that $\tilde\omega_{n,k,\Cal V}$ is
surjective. Then
\roster
\item
If $E_{n,k,\Cal V}^*$ is $(-2)$-regular, then
$
\reg{X}\leq(d-e+1)+\sum_{j=3}^{k}(j-2)\dim{V_j};
$
\item
If $E_{n,k,\Cal V}^*$ is $(-3)$-regular, then
$$
\reg{X}\leq (d-e+1)-2\dim{V_1}-\dim{V_2}+\sum_{j=4}^k(j-3)\dim{V_j}.
$$
\endroster
\endproclaim

\demo{Proof}
This is an easy consequence of Lemma 2.1,\,(2) (cf\. \cite{K2}
or \cite{L}).
\qed
\enddemo

\proclaim{Proposition 2.4}
Let $X$ be a $5$-dimensional smooth variety of degree $d$ in
$\Bbb P^8$.
\roster
\item
If $H^0(\Cal I_{X/\Bbb P^8}(2))\neq 0$, then $\reg{X}\le d+4$\rom;
\item
If $H^0(\Cal I_{X/\Bbb P^8}(2))=0$, then $\reg{X}\le d-5<d-e+1=d-2$.
\endroster
\endproclaim

\demo{Proof} We first deal with the case when $X$ is contained in a
hyperquadric $Q$. Note that all trisecant lines of $X$ are contained
in $Q$. As in \cite {L} and \cite{K2}, we consider the projection
$\pi_{\ell}\:X\to Y= \pi_{\ell}(X)\subset\Bbb P^6$, where
$\ell\subset\Bbb P^8$ is a generic line. Let $\ell\cap Q=\{p_1,p_2\}$.
If some fiber has three collinear points, then the line
through these three points contains $p_1$ or $p_2$. Let $\ell=\Bbb
P(V)$, $V=\Bbb CT_1\oplus\Bbb CT_2$, where $T_1$ and $T_2$ are linear
forms nonvanishing at $p_1$ and $p_2$.

Put $Y_k=\{y\in Y\mid\length{\pi_{\ell}^{-1}(y)}\geq k\}$. By Theorem
1.6, $\dim{Y_k}\leq6-k$ and $Y_6\subseteq Y_5\subseteq Y_4\subseteq
Y_3\subseteq Y_2\subset Y_1=Y=\pi_{\ell}(X)\subset\Bbb P^6$.
We will show that the morphism
$$
\bigoplus_{i=3}^5T_1^i\otimes\Cal O_{\Bbb P^6}(-i)
\oplus S^2(V)\otimes\Cal  O_{\Bbb P^6}(-2)\oplus V\otimes
\Cal O_{\Bbb P^6}(-1)\oplus\Cal O_{\Bbb P^6}\to\pi_{\ell\,*}\Cal O_X
\tag2.5
$$
is surjective.

We start with the following elementary lemma. From this lemma, we can choose
specific polynomials of degree $(n-1)$ to show $(n-1)$-normality of fibers of
generic projection from a given variety to a hypersurface.
\proclaim{Lemma 2.6}
Let $U,T_1,T_2$ be homogeneous coordinates on $\Bbb P^2$, and consider a collection
of $n+3$ points $p_i$,\,\,$i=1,\dots,n+3,\, n\ge3$.
\roster
\item
$p_i=(u_i,a,b)$, $1\le i\le n+1$ and the remaining two points $p_{n+2}, p_{n+3}$
are not contained in the line $aT_2-bT_1=0$.
\item
$p_i=(u_i,a,b)$, $1\le i\le n$ and the remaining three points $p_{n+1}, p_{n+2}, p_{n+3}$
are not contained in the line $aT_2-bT_1=0$.
\endroster
Suppose that $a\neq0$, $b\neq0$, $u_i\neq0$,
$i=1,\dots,n+1$. Then the points $p_i$, $i=1,\dots,n+3$ can be separated
by $(n+4)$ number of the monomials of degree $n$
$$
U^{n-j}{T_1}^{j},0\le j\le n,\ U^{n-1}T_2,\,U^{n-2}{T_2}^2,\,U^{n-2}T_1T_2,
$$
\endproclaim
\demo{Proof} We argue as in the proof of Lemma~3.3 in \cite{K2}.
The proof of cases (1) and (2) are almost same. So, we give a proof of (2) here.
By symmetry, it suffices to construct a form of degree $n$ vanishing
at all points except $p_n$ on the line and to construct a form
of degree $n$ vanishing at all points but $p_{n+3}$ off the line.

Consider the following system of polynomials of degree $n$:
$$
U^{n}+U^{n-1}(a_1T_1+a_2T_2)+U^{n-2}(a_3T_1^2+a_4T_1T_2+
a_5T_2^2)+\sum_{j=1}^{n-2}c_jU^{n-2-j}T_1^{2+j}.
$$
We observe that this form doesn't vanish identically on the line
$aT_2-bT_1=0$ containing the $n$ aligned points
$p_1,\dots,p_n$ because it doesn't vanish at the point $(1,0,0)$ on the line.
Substituting $T_1=a$, $T_2=b$, we get a system of polynomials
of degree $n-1$ in $U$:
$$
U^{n}+(a_{1}a+a_{2}b)U^{n-1}+(a_{3}a^2+a_{4}ab+a_{5}b^2)U^{n-2}
+\sum_{j=1}^{n-2} c_{j}a^{2+j}U^{n-2-j}.
$$
In order that our polynomial vanish at $p_1$, $p_2$,..., $p_{n-1}$, $p_{n}$,
it should be equal to $(U-u_1)(U-u_2)\cdots(U-u_{n-1})(U-u_n)$. Thus we get a
system of $n$ linear equations in the $n+3$ unknowns. This system
has a $3$-dimensional family of solutions, which allows us to pick a
solution passing through $p_{n+1}$, $p_{n+2}$.
Thus we constructed a form of degree $n$ passing though all points
except $p_{n+3}$. Therefore the point $p_{n+3}$ can be separated from the
other $n+2$ points. Similarly, we can construct a form of degree $n$ vanishing
at all points except $p_n$ on the line.
\qed
\enddemo

We return to the proof of Proposition~2.4. By Proposition~1.2, the
morphism
$$
T_1^3\otimes\Cal O_{\Bbb P^6}(-3)\oplus S^2(V)\otimes
\Cal O_{\Bbb P^6}(-2)\oplus V\otimes\Cal O_{\Bbb P^6}(-1)\oplus
\Cal O_{\Bbb P^6}\to\pi_{\ell\,*}\Cal O_X
$$
is surjective over the complement of the subvariety $Y_5$.

Let $y\in Y_5\setminus Y_6$. Using  Lemma 1.3 and
inequality~(1.0), one can show that the fiber $X_y$ is of
one of the following types:
\roster
\item
$X_y$ consists of five distinct points.
\itemitem{(i)}
$X_y$ consists of five distinct collinear points. In this case $X_y$
is $4$-normal and the monomials $T_1^iU^j$, $i+j=4$, $i,j\ge0$
generate a complete linear system of quartics;
\itemitem{(ii)}
$X_y$ consists of five distinct points only four of which are
collinear. In this case $X_y$ is $3$-normal, and to distinguish
between the points it suffices to use the monomials from
Lemma~2.6.(1) for $n=3$;
\itemitem{(iii)}
$X_y$ consists of five distinct points no four of which are collinear.
In this case $X_y$ is $2$-normal;
\item
Suppose that $X_y$ has multiple points and
$\Supp({X_y})$ consists of four distinct points $x_1,x_2,x_3,x_4$.\,\,(In this case
$X_y$ has only one double point in view of the inequality~(1.0))
\itemitem{(i)}
$\dim{\langle X_y\rangle}=1$. In this case $X_y$ is
$4$-normal by the monomials $T_1^iU^j$, $i+j=4$, $i,j\ge0$;
\itemitem{(ii)}
$\dim{\langle X_y\rangle}=2$. In this case $X_y$ is
$2$-normal except $X_y$ has a collinear subscheme of length $4$. In this case,
it is $3$-normal by using the
monomials from Lemma~2.6.(1) for $n=3$;
\item
$\Supp({X_y})$ consists of at most three distinct points. This
case is impossible in view of the inequality~(1.0).
\endroster

Thus the morphism
$$
\bigoplus_{i=3}^4\bigl(T_1^i\otimes\Cal O_{\Bbb P^6}(-i)\bigr)\oplus
S^2(V)\otimes\Cal O_{\Bbb P^6}(-2)\oplus V\otimes\Cal O_{\Bbb P^6}(-1)
\oplus\Cal O_{\Bbb P^6}\to\pi_{\ell\,*}\Cal O_X
$$
is surjective over all points $y\notin Y_6$.

By Theorem~1.6, for $y \in Y_6$, a finite set, the fiber $X_y$ consists
of six distinct points. Furthermore, using Proposition 1.2 and Lemma 2.6  we know
that the vector space of monomials
$$
U^{5-j}{T_1}^{j},0\le j\le 5,\ U^{4}T_2,\,U^{3}{T_2}^2,\,U^{3}T_1T_2
$$
separate points in $X_y$ for all $y \in Y_6$.

(Note that $U$ can be chosen as a linear
form in $\Bbb P^6$ not through $y \in Y_6$ such that $U$,$T_1$, $T_2$ play a role as
coordinates in ${\langle \ell,y\rangle}\simeq \Bbb P^2$). Thus the morphism~(2.5) is
surjective over $\Bbb P^6$ and, by Lemma~2.3\,(1), $\reg{X}\le \deg{X}+4$.

For~(2), suppose that $X$ is not contained in a
hyperquadric $Q$. Note that $X$ is linearly normal (Zak's Theorem) and
$H^1(\Cal O_X)=0$ (Barth's Theorem). On the other hand,
Proposition~1.2  and Corollary 1.3 show that the natural morphism
$$
S^3(V)\otimes\Cal O_{\Bbb P^6}(-3)\oplus S^2(V)\otimes
\Cal O_{\Bbb P^6}(-2)\oplus V\otimes\Cal O_{\Bbb P^6}(-1)\oplus
\Cal O_{\Bbb P^6}\to\pi_{\ell\,*}\Cal O_X
$$
is surjective over all $y\in Y$ for which the fiber $X_y$ does not
contain five collinear points. By Theorem~1.6,
$\dim{\{y\in Y_5\mid X_y\text{\ has five aligned points}\}}\le1$.
If $\deg{(X_y)}=5$ and $\dim{\langle X_y\rangle}=1$,
then $\dim{\{q\in \ell\mid \langle q,y\rangle\text{\ is a five secant line of}\, X_y \}}$
may be equal to one and thus, we need  quartic polynomials either $\sum_{i=0}^4a_iU^{4-i}T_1^i$ or
$\sum_{i=0}^4a_iU^{4-i}T_2^i$, where $U$ is a linear form in $\Bbb P^6$ which
does not vanish at a point $y \in Y_5\setminus Y_6$. Thus, letting $V_4$ be $\{T_1^4,T_2^4\}$ as a
subspace of $S^4(V)$,
$$
V_4\otimes
\Cal O_{\Bbb P^6}(-4)\bigoplus_{i=0}^3S^i(V)\otimes
\Cal O_{\Bbb P^6}(-i)\to\pi_{\ell\,*}\Cal O_X
$$

Consider now the fibers over the finite set $Y_6$. As in the proof of (1),
we know
that the vector space of monomials
$$
U^{5-j}{T_1}^{j},0\le j\le 5,\ U^{4}T_2,\,U^{3}{T_2}^2,\,U^{3}T_1T_2
$$
separate distinct six points in $X_y$ for all $y \in Y_6$. Therefore the morphism
$$
T_1^5\otimes\Cal O_{\Bbb P^6}(-5)\oplus
V_4\otimes\Cal O_{\Bbb P^6}(-4)
\bigoplus_{i=0}^3S^i(V)\otimes
\Cal O_{\Bbb P^6}(-i)\to\pi_{\ell\,*}\Cal O_X
$$
is surjective over $\Bbb P^6$.
By Lemma~2.1\,(2), the dual of the kernel $E_{5,6,\Cal V}$ is
$(-3)$-regular. We get $$\reg{X}\le\deg{X}-5 < d-e+1=d-2$$ from Lemma~2.3\,(2).
\qed
\enddemo

\proclaim{Corollary 2.7} Let $X^5\subset\Bbb P^{8}$, be a smooth projective variety.
If $X$ has $(d-4)$-secant line
$l\not\subset X$, then $X$ \rom(and $l$\rom) is contained in a quadric
hypersurface $Q\subset\Bbb P^{8}$.
\endproclaim

\demo{Proof} Since $m$-regularity of $X$ implies that the degrees of
defining equations are bounded by $m$ (cf\.~\cite{Mu1, Lecture~14}),
this follows immediately from Remark~2.7\,(2).
\enddemo

\proclaim{Proposition 2.8}
Let $X$ be a smooth projective variety of dimension $6$ in $\Bbb P^9$.
\roster
\item
If $H^0(\Cal I_{X/\Bbb P^9}(2))\neq 0$, then $\reg{X}\le d+8$\rom;
\item
If $H^0(\Cal I_{X/\Bbb P^9}(2))=0$, then $\reg{X}\le d$.
\endroster
\endproclaim
\demo{Proof}
We consider the projection $\pi_{\ell}\:X\to Y= \pi_{\ell}(X)\subset\Bbb P^7$, where
$\ell\subset\Bbb P^9$ is a generic line, where $\ell=\Bbb P(V)$, $V=\Bbb CT_1\oplus\Bbb CT_2$.
For a proof of (1), the arguments used to prove Proposition~2.6.(1) can show that
the morphism
$$
\bigoplus_{i=3}^6T_1^i\otimes\Cal O_{\Bbb P^{7}}(-i)
\bigoplus_{i=0}^2S^i(V)\otimes\Cal O_{\Bbb P^{7}}(-i)\to
\pi_{\ell\,*}\Cal O_X
$$
is surjective over $\Bbb P^7$. Indeed, for any reduced fiber $X_y$, we have $\deg{(X_y)}\le7$ and
as in the proof of Proposition 2.6, each point of $X_y$ can be separated by
the vector space of monomials
$$
U^{6-j}{T_1}^{j},0\le j\le 6,\ U^{5}T_2,\,U^{4}{T_2}^2,\,U^{4}T_1T_2.
\tag*
$$
Now, it is enough to consider nonreduced fiber $X_y$ such that $\dim{\langle X_y\rangle}=2$
for $y \in Y_5$. Now, suppose that $\dim{\langle X_y\rangle}=2$ and $\deg{(X_y)}=5$.
Then, $X_y$ is $2$-normal except $X_y$ has a $4$-secant line (by Corollary 1.3. (a)),
in which case $X_y$ is $3$-normal by the same arguments as those in Proposition 2.6.

Finally, suppose $\dim{\langle X_y\rangle}=2$, $\deg{(X_y)}=6$ and
$X_y$ has a nonreduced point.
\roster
\item
Suppose $X_y$ has $5$-secant line.
In this case, $X_y$ is $4$-normal by monomials
$U^{4-j}{T_1}^{j},0\le j\le 4,\ U^{3}T_2,\,U^{3}{T_2}^2,\,U^{2}T_1T_2$.
\item
Suppose $X_y$ has no $5$-secant line. Then it is $3$-normal. In particular,
\itemitem{(i)}
$\Supp({X_y})$ consists of five distinct points $x_1,x_2,x_3,x_4,x_5$.\,\,(In this case
$X_y$ has only one double point in view of the inequality~(1.0)) This case is
$3$-normal by monomials
$U^{3-j}{T_1}^{j},0\le j\le 3,\ U^{2}T_2,\,U{T_2}^2,\,UT_1T_2$.
\itemitem{(ii)}
$\Supp{(X_y)}$ consist of four points $x_1,x_2,x_3,x_4$. Note that the Zariski tangent
space of $X_y$ at a nonreduced point might be two-dimensional under the generic projection
$\pi_{\Lambda}\:X^6\to Y^6\subset\Bbb P^7$ (see Lemma 1.5). However, in this case it is
impossible in view of the inequality~(1.0) and (see page 190 \cite {Ma1}).
So, $\Supp{(X_y)}$ has no four points.
\itemitem{(iii)}
$\Supp({X_y})$ consists of at most three distinct points. This
case is also impossible in view of the inequality~(1.0).
\endroster
Therefore, the morphism
$$
\bigoplus_{i=3}^6T_1^i\otimes\Cal O_{\Bbb P^{7}}(-i)
\bigoplus_{i=0}^2S^i(V)\otimes\Cal O_{\Bbb P^{7}}(-i)\to
\pi_{\ell\,*}\Cal O_X
$$
is surjective and $\reg{X}\le d-2+(1+2+3+4)=d+8$.

For a proof of (2), $X_y$ is $3$-normal
for $y \notin Y_5$, and  for $\deg{X_y}=5$ $X_y$ is $4$-normal and fails to be $3$-normal iff
it is contained in a line. However, since $\{q\in \ell\mid \langle q,y\rangle\text
{\ is a five secant line of}\,\, X_y \}$ might be {\it $two$}-dimensional
it is clear that letting $V_4$ be $\{T_1^4,T_2^4\}$,
$$
V_4\otimes
\Cal O_{\Bbb P^6}(-4)\bigoplus_{i=0}^3S^i(V)\otimes
\Cal O_{\Bbb P^6}(-i)\to\pi_{\ell\,*}\Cal O_X
$$
is surjective over the complement of $Y_6$.

Next, assume that $\dim{\langle X_y\rangle}=1$ for $y\in Y_6\setminus Y_7$.
Then it is $5$-normal and  it is enough to choose $V_5=\{T_1^5,T_2^5\}$ as a subspace of $S^5(V)$
(since $\dim{Y_6}\le 1$).
Finally, suppose either $\deg{(X_y)}=6$, $\dim{\langle X_y\rangle}=2$ or $\deg{(X_y)}=7$.
it is reduced to the same arguments as those used in a proof of (1).
As a consequence,
$$
T_1^6\otimes\Cal O_{\Bbb P^7}(-6)\oplus
V_5\otimes\Cal O_{\Bbb P^7}(-5)\oplus
V_4\otimes\Cal O_{\Bbb P^7}(-4)
\bigoplus_{i=0}^3S^i(V)\otimes
\Cal O_{\Bbb P^7}(-i)\to\pi_{\ell\,*}\Cal O_X
$$
is surjective over $\Bbb P^7$.
Since $H^0(\Cal I_{X/{\Bbb P}}(2))=H^1(\Cal O_X)=0$ and $X$ is linearly normal,
by Lemma 2.3.(2), $\reg{X}\le d$.
\enddemo

\remark{\bf Remark 2.9}
For a smooth variety of dimension $n$ in $\Bbb P^{n+3},\,\, n\ge7$, Hartshorne's conjecture
states that $X$ should be a complete intersection. So, if $X$ is a complete intersection
whose homogenious ideal is generated by three polynomials of degrees $d_1$, $d_2$ and $d_3$
then, by a standard computation  $\reg{X}=d_1+d_2+d_3-2$.
\endremark

\proclaim{Theorem 2.10} Let $X\subset\Bbb P^N$ be a smooth projective
variety of dimension~$5$ and codimension $e\ge4$. Then
$\reg{X}\le(d-e+1)+10$.
\endproclaim

\demo{Proof} As in Proposition 2.4, we consider
the projection $\pi_{\Lambda}\:X\to Y\subset\Bbb P^6$ with center at a
generic linear subspace $\Lambda=\Bbb P(V)\simeq \Bbb P^{e-2}$, $V=\Bbb C\cdot T_7\oplus
\Bbb C\cdot T_8\oplus\dots\oplus\Bbb C\cdot T_N$, $\Lambda\cap X=\varnothing$.
Put $Y_k=\{y\in Y\mid\length{\pi_{\Lambda}^{-1}(y)}\geq k\}$ as before.
The morphism
$\bigoplus_{i=0}^2S^i(V)\otimes\Cal O_{\Bbb P^6}(-i)\to
\pi_{\Lambda\,*}\Cal O_X$ is surjective for all $y\notin Y_4$ because
the fibers $X_y$ are $2$-normal.
Let's consider all cases according to $\dim{\langle X_y\rangle}$.
\roster
\item
If $\dim{\langle X_y\rangle}\ge3$, then  $X_y$ is $2$-normal except that $X_y$
consists of distinct six points with $4$-secant line because ${\deg(X_y)}\le6$.
Such a special fiber $X_y$ can be brought
into a general position by cubic forms
$U^{3-i}L_1^i,\, 0\le i\le3,\, UT_{ij}, \,T_{ij}\in S^2(V)$,
where $U$ is a linear form in $\Bbb P^6$ which does not vanish at $y$ and
$L_1$ is a linear form in $\Lambda$ nonvanishing on the $4$-secant line.
\item
If $\dim{\langle X_y\rangle}\le2$, then the fibers $X_y$ are well described
in Proposition 2.4. However, $\{q\in \Lambda\mid \langle q,y\rangle\text
{\ is a four secant line of}\,\, X_y \}$ might be {\it $two$}-dimensional and $\{q\in \Lambda\mid \langle q,y\rangle\text
{\ is a five secant line of}\,\, X_y \}$ might be {\it $one$}-dimensional.
Thus, it suffices to choose three linear forms $L_1,L_2$ and $L_3$ in $\Lambda$
such that $V_3=\{L_1^3, L_2^3, L_3^3\}\subset S^3(V)$, $V_4=\{L_1^4, L_2^4\}\subset S^4(V)$,
$V_5=\{L_1^5\}\subset S^5(V)$.
\endroster
Summing up, we get a surjective morphism
$$
V_5\otimes\Cal O_{\Bbb P^6}(-5)\oplus
V_4\otimes\Cal O_{\Bbb P^6}(-4)\oplus
V_3\otimes\Cal O_{\Bbb P^6}(-3)
\bigoplus_{i=0}^2S^i(V)\otimes
\Cal O_{\Bbb P^6}(-i)\to\pi_{\Lambda\,*}\Cal O_X
\tag**
$$
over $\Bbb P^6$. Therefore, by Lemma 2.3\,(1),
$\reg{X}\le\bigl(\deg{X}-\codim{X}+1\bigr)+10$.
\qed
\enddemo

\proclaim{Corollary 2.11} Let $X\subset\Bbb P^N$ be a smooth projective
variety of dimension~$5$ and codimension $e\ge4$. If $X$ is linearly normal and
$H^0(\Cal I_X(2))=H^1(\Cal O_X)=0$, then
$\reg{X}\le(d-e+1)-2(e-1)-\frac{e(e-1)}{2}+4$.
\endproclaim
\demo{Proof}
It is easily computed from (**), Lemma 2.1.(2) and Lemma 2.3.(2).
\enddemo
\remark{\bf Remark~2.12}
Let $X$ be a smooth projective variety of codimension $e$, and $I_X$ be
the saturated ideal of $X$. Suppose  $\dim{X}\le5$ and $(e-1)$ defining
equations out of minimal generators of $I_X$ have degree two. Then for all $y\in Y$,
$X_y$ has no trisecant line and $\deg{X_y}$ is at most $6$.
This implies that $X_y$ is $2$-normal and the following morphism
$
\bigoplus_{i=0}^2S^i(V)\otimes
\Cal O_{\Bbb P}(-i)\to\pi_{\Lambda\,*}\Cal O_X
$
is surjective for all $y\in Y$.
By Lemma 2.3.(1), $\reg{X}\le d-e+1$.
\endremark

\proclaim{Theorem 2.13} Let $X\subset\Bbb P^N$ be a smooth projective
variety of dimension~$6$ and codimension $e\ge4$. Then
$\reg{X}\le(d-e+1)+20$.
\endproclaim

\demo{Proof} As in Proposition 2.8, we consider all possible fibers according to
$\deg{X_y}$ and $\dim{\langle X_y \rangle}$. Consider the projection
$\pi_{\Lambda}\:X\to Y\subset\Bbb P^7$
with center at a generic linear subspace $\Lambda=\Bbb P(V)\simeq \Bbb P^{e-2}$,
$V=\Bbb C\cdot T_8\oplus
\Bbb C\cdot T_8\oplus\dots\oplus\Bbb C\cdot T_N$, $\Lambda\cap X=\varnothing$.
Note that $\deg{X_y}\le 7$ for all $y\in Y$. By Proposition 1.2,
the fibers $X_y$ are $2$-normal if either $y\in Y\setminus Y_4$ or
$\dim{\langle X_y \rangle}\ge5$. In addition, it is easy to consider the normality of the fibres $X_y$ which consists of distinct seven points for $y\in Y_7$ as before.
Let's consider all the remaining cases according to $\dim{\langle X_y\rangle}$
and $\deg{(X_y)}\le 6$.
\roster
\item
If $\dim{\langle X_y\rangle}=4$ then it is $2$-normal.
\item
Suppose $\dim{\langle X_y\rangle}=3$.
If $X_y$ is a fiber of degree $6$ with 4-secant line, it is $X_y$ is $3$-normal by cubic polynomials
$U^{3-i}L_1^i,\, 0\le i\le3,\, UT_{ij}, \,T_{ij}\in S^2(V)$,
where $U$ is a linear form in $\Bbb P^7$ which does not vanish at $y$ and
$L_1$ is a linear form in $\Lambda$ nonvanishing on the $4$-secant line.)
Otherwise, $X_y$ is $2$-normal.
\item
If $\dim{\langle X_y\rangle}\le2$, then the fibers $X_y$ are well described
in Proposition 2.4. However, as in Proposition 2.10, $\{q\in \Lambda\mid \langle q,y\rangle\text
{\ is a four secant line of}\,\, X_y \}$ might be {\it $three$}-dimensional, $\{q\in \Lambda\mid \langle q,y\rangle\text{\ is a five secant line of}\,\, X_y \}$ might be {\it $two$}-dimensional and
$\{q\in \Lambda\mid \langle q,y\rangle\text{\ is a six secant line of}\,\, X_y \}$ might be {\it $one$}-dimensional.
Thus, it suffices to choose four linear forms $L_1,L_2$ and $L_3$ $L_4$ in $\Lambda$
such that $V_3=\{L_1^3, L_2^3, L_3^3, L_4^3\}\subset S^3(V)$, $V_4=\{L_1^4, L_2^4, L_3^4\}\subset S^4(V)$,
$V_5=\{L_1^5, L_2^5\}\subset S^5(V)$, $V_6=\{L_1^6\}\subset S^6(V)$   .
\endroster
Finally, we get a surjective morphism
$$
\bigoplus_{i=3}^6V_i\otimes
\Cal O_{\Bbb P^7}(-i)
\bigoplus_{i=0}^2S^i(V)\otimes
\Cal O_{\Bbb P^7}(-i)\to\pi_{\Lambda\,*}\Cal O_X
\tag***
$$
over $\Bbb P^7$. Therefore, by Lemma 2.3\,(1),
$\reg{X}\le\bigl(\deg{X}-\codim{X}+1\bigr)+20$.
\qed
\enddemo

\proclaim{Corollary 2.14} Let $X\subset\Bbb P^N$ be a smooth projective
variety of dimension~$6$ and codimension $e\ge4$. If $X$ is linearly normal and
$H^0(\Cal I_X(2))=H^1(\Cal O_X)=0$, then
$\reg{X}\le(d-e+1)-2(e-1)-\frac{e(e-1)}{2}+10$.
\endproclaim
\demo{Proof}
It is easily computed from (***), Lemma 2.1.(2) and Lemma 2.3.(2).
\enddemo

\head
\bf \S3. Castelnuovo regularity for integral projective surfaces
\endhead

The known results in and approaches to the Castelnuovo regularity
problem for integral projective varieties are quite different from
those in the smooth cases.

For example, for an arbitrary toric variety
$X\subset\Bbb P^N$ the best known bound is $\reg{X}\le(N+1)\cdot
\deg{X}\cdot\codim{X}$, but if $X$ is a toric variety of codimension
two there is a much better bound $\reg{X}\le\deg{X}$ (cf\.~\cite{PS}).
Another bound on regularity is known in the case of Buchsbaum
varieties, (cf\.~\cite{HM} and \cite{SV}).
                             
In this section we try to extend techniques used in section \S2 to
integral projective varieties. Note that some arguments in section \S2
can not be applied to integral projective varieties because they depend on
the Kodaira vanishing theorem and information on generic linear
projections of smooth varieties to a hypersurface in projective space. However for
integral curves, the sharp bound $\reg{X}\le d-e+1$ is proved and
classification of extremal curves for which $\reg{X}=d-e+1$ is
given in~\cite{GLP}. Thus the case of integral projective surfaces
is the simplest case when the regularity conjecture is still open. By
Lemma~1.1, for any integral projective variety $X$ one has
$$
\reg{(X\cap H)}\le\reg{X}\le\reg{(X\cap H)}+h^1\bigl(\Cal I_X(\reg{(X\cap H)}-1)\bigr),
$$
where $H$ is a generic hyperplane in $\Bbb P^N$. Hence, for integral
projective surfaces, $m$-normality implies $(m+1)$-regularity for
$m\ge d-e$.
In addition, for a generic hyperplain section $Y= X\cap H$,
since $Y$ is $(d-e+1)$-regular, so by Lemma 1.1 (b), $X$ is $(d-e+1)+h^1({\Cal I}_{X}(d-e))$-regular.
On the other hand,
$$
h^1({\Cal I}_{X}(d-e))\le h^0({\Cal O}_{\Bbb P^N}(d-e))-\chi({\Cal I}_{X}(d-e))
$$
But $\chi({\Cal I}_{X}(n))=\chi({\Cal O}_{\Bbb P^N}(n))-\chi({\Cal O}_{X}(n))$
for all $n$, and $\chi({\Cal O}_{\Bbb P^N}(n))=h^0({\Cal O}_{\Bbb P^N}(n))$ for all $n\ge 0$.
Hence, $h^1({\Cal I}_{X}(d-e))\le \chi({\Cal O}_{X}(d-e))=\Cal P_X(d-e)$ and
$$
\reg{X}\le (d-e+1)+h^1({\Cal I}_{X}(d-e))\le (d-e+1)+ \Cal P_X(d-e)
$$
which is a cubic polynomials in $d$.
In this section, we give a quadratic bound in $d$ for $\reg{X}$ for arbitrary
integral projective surface $X$.

We proceed with using methods from \S2 to recover the construction
in~\cite{G}, to justify Greenberg's unsubstantiated claims, and to
improve his regularity bound in the case of integral projective
surfaces.

Greenberg's main idea consists in considering a general linear
projection from $X$ to $\Bbb P^2$ and applying the Eagon-Northcott
complex to an exact sequence of vector bundles.

Let $X$ be a nondegenerate integral complex projective surface of
degree $d$ in $\Bbb P^N$ (we do not assume $X$ to be locally
Cohen-Macaulay). Let $\Lambda^{N-3}\subset\Bbb P^N$ be a generic
linear subspace, $\Lambda\cap X=\varnothing$, and let
$p_1\:Bl_{\Lambda}\Bbb P^N\to\Bbb P^N$ be the blowing up of $\Lambda$.
Let $p_2\:Bl_{\Lambda}\Bbb P^N\to\Bbb P^2$, $p_2=\pi_{\Lambda}\circ
p_1$, where $\pi_{\Lambda}\:\Bbb P^N\dasharrow\Bbb P^2$ is the
projection with center at $\Lambda$.  Without loss of generality, we
can choose homogeneous coordinates $T_0,T_1,\dots,T_N$ in $\Bbb P^N$
so that $\Lambda=Z(T_0,T_1,T_2)$ is defined by vanishing of $T_0$,
$T_1$ and $T_2$. Putting $V=\Bbb C\cdot T_3\oplus\Bbb C\cdot T_4\oplus
\cdots\oplus\Bbb C\cdot T_N$, and we also have $\Lambda=\Bbb P(V)$.

Note that $Bl_{\Lambda}\Bbb P^N=\Bbb P\bigl(\Cal O_{\Bbb P^2}(1)\oplus
(V\otimes\Cal O_{\Bbb P^2})\bigr)=\{(x,q)\mid x\in L_q=
\langle\Lambda,q\rangle,\ q\in\Bbb P^2\}$.

As in \cite{K2} or \cite{L}, consider the diagram
$$
\CD
Bl_{\Lambda}\Bbb P^N=\Bbb P\bigl(\Cal O_{\Bbb P^2}(1)\oplus(V\otimes
\Cal O_{\Bbb P^2})\bigr)@>p_2>>\Bbb P^2\\
@VVp_1V@.\\
X\subset\Bbb P^N@.
\endCD
$$
From the choice of $\Lambda$ and the definition of degree, it follows that
$\pi_{\Lambda}\:X\to\Bbb P^2$ is a $d:1$ morphism, i.e\. the
fibers $\pi^{-1}_{\Lambda}(q)$ have length $d$ for all $q\in\Bbb P^2$.

Consider the morphism
$$
p_{2\,*}\bigl(p_1^{\ast}\Cal O_{{\Bbb P}^N}(k)\bigr)
\overset{\omega_{2,k}}\to\longrightarrow
p_{2\,*}\bigl(p_1^*\Cal O_X(k)\bigr)
$$
Note that $p_{2\,*}\bigl(p_1^*\Cal O_X(1)\bigr)=\pi_{\Lambda\,*}\Cal O_X(1)$,
$\pi_{\Lambda\,*}\Cal O_X(k)
=\pi_{\Lambda\,*}\Cal O_X\otimes\Cal O_{\Bbb P^2}(k)$
and $p_1^*\Cal O_{\Bbb P^N}(1)
=\Cal O_{\Bbb P(\Cal E)}(1)$ (the tautological line  bundle),
where $\Cal E=\Cal O_{\Bbb P^2}(1)\oplus(V\otimes\Cal O_{\Bbb P^2})$.

The main issue is to prove the surjectivity of $\omega_{2,k}$ for a
suitable $k>0$. By Nakayama's lemma, it suffices to show that for all
$q\in\Bbb P^2$ the upper arrow in the commutative diagram
$$
\CD
p_{2\,*}\bigl(p_1^*\Cal O_{\Bbb P^N}(k)\bigr)\otimes\Bbb C(q)@>
\omega_{2,k}\otimes\Bbb C(q)>>p_{2\,*}\bigl(p_1^*\Cal O_X(k)\bigr)
\otimes\Bbb C(q) \\
@V\cong VV@V\cong VV\\
H^0\bigl(L_q,\Cal O_{L_q}(k)\bigr)@>>>H^0\bigl(L_q,
\Cal O_{\pi^{-1}_{\Lambda}(q)}(k)\bigr)
\endCD
$$
is surjective for some $k>0$. Equivalently, using the bottom arrow, it is
enough to show that the finite scheme $\pi^{-1}_{\Lambda}(q)$ of
length $d$ in $L_q=\langle\Lambda,q\rangle$ is $k$-normal for a
suitable number $k>0$. We need the following lemma.

\proclaim{Lemma 3.1} Let $X$ be an integral projective variety
of dimension $n$ and degree $d$ in $\Bbb P^N$. Let $r\le N-n$ be a
natural number, and let $\Bbb P^r\subset\Bbb P^N$ be a linear
subspace. We put $Y=\Bbb P^r\cap X$. Then
$$
\deg{Y}\le d-(N-n-r)
$$
\endproclaim
\demo{Proof} Adding $(N-n-r)$ general points of $X$ to $Y$, we see that
the lemma follows from the generalized Bezout theorem.\qed
\enddemo
\bigskip
From Lemma~3.1, we get $\langle\pi^{-1}_{\Lambda}(q)\rangle=\Bbb P^{N-2}$ and
$\pi^{-1}_{\Lambda}(q)$ is $k$-normal for all $q\in\Bbb P^2$ and all $k\ge d-(N-2)$ by
Proposition~1.2.

Recall that
$$
\multline
p_{2\,*}\bigl(p_1^*\Cal O_{\Bbb P^N}(k)\bigr)
=Sym^k\bigl(\Cal O_{\Bbb P^2}(1)\oplus(V\otimes\Cal O_{\Bbb P^2})
\bigr)\\
=\Cal O_{\Bbb P^2}(k)\oplus V\otimes\Cal O_{\Bbb P^2}(k-1)\oplus S^2(V)
\otimes\Cal O_{\Bbb P^2}(k-2)\oplus\cdots\oplus S^k(V)\otimes
\Cal O_{\Bbb P^2},
\endmultline
$$
where $S^i(V)$ is the $i$-th symmetric power of $V$. After twisting by
$(-k)$, we get an exact sequence
$$
0\to E_{2,k}\to S^k(V)\otimes\Cal O_{\Bbb P^2}(-k)\oplus\cdots
\oplus V\otimes\Cal O_{\Bbb P^2}(-1)\oplus\Cal O_{\Bbb P^2}
\overset{\omega_{2,k}}\to\longrightarrow\pi_{\Lambda\,*}\Cal O_X\to0,
\tag3.2
$$
where  $E_{2,k}=\Ker{\omega_{2,k}}$. Note that $\pi_{\Lambda\,*}\Cal
O_X$ is a vector bundle of rank $d$ over $\Bbb P^2$ because the
projection $\pi_{\Lambda}\:X\to\Bbb P^2$ is flat (all fibers of
$\pi_{\Lambda}$ have the same length $d$).

Now, following Greenberg (\cite{G}), we generalize this construction.

\proclaim{Theorem 3.3} Let $X$ be an integral projective surface of
degree $d$ and codimension $e=N-2$ in $\Bbb P^N$. Then
$\reg{X}\le(d-e+1)d-(2e+1)$.
\endproclaim
\demo{Proof}
Let $\Cal F$ be a vector bundle $S^k(V)\otimes\Cal O_{\Bbb P^2}(-k)
\oplus\cdots\oplus V\otimes\Cal O_{\Bbb P^2}(-1)\oplus\Cal O_{\Bbb P^2}$
of rank $f$. The exact sequence $0\to E_{2,k}
\to\Cal F\overset{\omega_{2,k}}\to\longrightarrow\pi_{\Lambda\,*}
\Cal O_X\to0$ gives rise to an Eagon-Northcott complex (see page 494 \cite{GLP})
$$
\multline
0\to\bigwedge^f\Cal F\otimes
S^{f-d-1}(\pi_{\Lambda\,*}\Cal O_X)^{\lor}
\otimes\det{(\pi_{\Lambda\,*}\Cal O_X )^{\lor}}
\to\cdots\\
\cdots\to\bigwedge^{d+i+1}\Cal F@!@!\otimes@!@! S^i(\pi_{\Lambda\,*}
\Cal O_X)^{\lor}\otimes\det{(\pi_{\Lambda\,*}\Cal O_X)^{\lor}}
\cdots\to\bigwedge^{d+2}\Cal F\otimes(\pi_{\Lambda\,*}
\Cal O_X)^{\lor}\otimes\det{(\pi_{\Lambda\,*}\Cal O_X )^{\lor}}\\
\overset{\varphi_1}\to\to\bigwedge^{d+1}\Cal F\otimes
\det{(\pi_{\Lambda\,*}\Cal O_X)^{\lor}}\overset{\varphi_0}\to\to
\Cal F\overset{\omega_{2,k}}\to\longrightarrow\pi_{\Lambda\,*}
\Cal O_X\to0
\endmultline
$$
Note that if $H^1\bigl(\Bbb P^2, E_{2,k}(m)\bigr)=0$ for some $m>0$, then
$X$ is $m$-normal from the sequence (3.2).
By chopping, we get the following two exact sequences of sheaves over
$\Bbb P^2$:
$$
\gathered
0\to\Ker{\varphi_0}\to\bigwedge^{d+1}\Cal F\otimes
\det{(\pi_{\Lambda_*}\Cal O_X} )^{\lor}\overset{\varphi_0}\to\to
E_{2,k}\to0;\\
0\to\Ker{\varphi_1}\to\bigwedge^{d+2}\Cal F
\otimes\det{(\pi_{\Lambda\,*}\Cal O_X)^{\lor}}\otimes
(\pi_{\Lambda\,*}\Cal O_X)^{\lor}\overset{\varphi_1}\to\to
\Ker{\varphi_0}\to0.
\endgathered
$$
Since $\bigwedge^{d+1}\Cal F\otimes
\det({\pi_{\Lambda\,*}\Cal O_X})^{\lor}$ is a direct sum of line bundles
of the form $\Cal O_{\Bbb P^2}(t)$, we have
$
h^1\bigl(\Bbb P^2,\wedge^{d+1}\Cal F\otimes
\det({\pi_{\Lambda\,*}\Cal O_X)^{\lor}(m)}\bigr)=0
$ for any $m\in  {\Bbb  Z}$. So, $h^2(\Bbb P^2, \Ker{\varphi_0(m)})=0$
implies $h^1\bigl(\Bbb P^2, E_{2,k}(m)\bigr)=0$. The remaining part is to
find out an integer $m_0$ such that $h^2(\Bbb P^2, \Ker{\varphi_0(m)})=0$
for all $m\ge m_0$.

By the way, from the second short exact sequence, it suffices
to find out an integer $m_0$ such that $h^2\bigl(\Bbb P^2,\bigwedge^{d+2}\Cal F
\otimes(\pi_{\Lambda\,*}\Cal O_X)^{\lor}\otimes\det{\bigl(\pi_{\Lambda\,*}\Cal O_X)^{\lor}}
(m)\bigr)=0$ for all $m\ge m_0$. On the other hand, one has
$$
\multline
h^2\bigl(\Bbb P^2,\wedge^{d+2}{\Cal F}\otimes(\pi_{\Lambda\,*}\Cal O_X)
^{\lor}\otimes\det{(\pi_{\Lambda\,*}\Cal O_X)^{\lor}}(m)\bigr)\\
=h^2\bigl(\Bbb P^2,\wedge^{d+2}\Cal F\otimes(\pi_{\Lambda\,*}\Cal O_X)
^{\lor}\otimes\Cal O_{\Bbb P^2}(-c_1+m)\bigr)\\
=h^0\bigl(\Bbb P^2,\wedge^{d+2}\Cal F^{\lor}\otimes\pi_{\Lambda\,*}
\Cal O_X\otimes\Cal O_{\Bbb P^2}(c_1-m-3)\bigr)
\endmultline
$$
where $c_1$ is
the first Chern class of $\pi_{\Lambda\,*}\Cal O_X$. Furthermore,
since $\pi_{\Lambda}$ is a finite affine morphism, we have
$$
H^0\bigl(\Bbb P^2,\wedge^{d+2}\Cal F^{\lor}\otimes\pi_{\Lambda\,*}
\Cal O_X\otimes\Cal O_{\Bbb P^2}(c_1-m-3)\bigr)
=H^0\bigl(X,\pi_{\Lambda}^*(\Lambda^{d+2}\Cal F^{\lor})
\otimes\Cal O_X(c_1-m-3)\bigr)
$$
and for $m\ge k(d+2)+c_1-2$, $h^0\bigl(X,\pi_{\Lambda}^*(\Lambda^{d+2}\Cal F^{\lor})
\otimes\Cal O_X(c_1-m-3)\bigr)=0
$
because
$\pi_{\Lambda}^*(\wedge^{d+2}\Cal F^{\lor})\otimes\Cal O_X(c_1-m-3)$
is a sum of line bundles $\Cal O_X(t)$ for $t<0$. To complete the proof
we need the following

\proclaim{Lemma 3.4}
$$
c_1(\pi_{\Lambda\,*}\Cal O_X)\le-d.
$$
\endproclaim
\demo{Proof} Since the first Chern class is stable with respect to
taking a general hyperplane section, we may assume that $X$ is an
integral curve of arithmetic genus $\rho_a$. Therefore, by the
Riemann-Roch theorem for vector bundles, $\chi(\pi_{\Lambda\,*}
\Cal O_X)=\rank{\pi_{\Lambda\,*}\Cal O_X}+c_1(\pi_{\Lambda_*}\Cal O_X)
=d+c_1(\pi_{\Lambda\,*}\Cal O_X)$. Since $\pi_{\Lambda}$ is a finite
morphism, $\chi(\pi_{\Lambda\,*}\Cal O_X)=\chi(\Cal O_X)=1-\rho_a$. Since
$X$ may be assumes to be singular, $\rho_a\ge1$ and
$c_1(\pi_{\Lambda\,*}\Cal O_X)=-d+1-\rho_a\le -d$.
\qed
\enddemo

Let's return to Theorem 3.3. From Lemma 3.4,
$H^1\bigl(\Bbb P^2, E_{2,k}(m)\bigr)=0$ for $m\ge k(d+2)-d-2$ and
$X$ is $m$-normal for $m\ge k(d+2)-d-2$.
As we already mentioned, in view of Proposition~1.2 and Lemma~3.1,
one can take $k=d-e\ge2$.  Then, Lemma~1.1\,(2) shows
that
$$
\reg{X}\le (d-e)(d+2)-d-1 = (d-e+1)d-(2e+1).
$$
\qed
\enddemo

\enddocument